\documentstyle[12pt]{article}
\newcommand{\sect}[1]{\section{#1}\setcounter{equation}{0}}

\font\mbn=msbm10 scaled \magstep1
\font\mbs=msbm7 scaled \magstep1
\font\mbss=msbm5 scaled \magstep1
\newfam\mbff
\textfont\mbff=\mbn
\scriptfont\mbff=\mbs
\scriptscriptfont\mbff=\mbss\def\mbf{\fam\mbff}
\def\Re{{\mbf R}}

\def\Z{{\mbf Z}}
\def\Co{{\mbf C}}
\def\To{{\mbf T}}

\def\H{{\mbf H}}
\newtheorem{Th}{Theorem}[section]

\newtheorem{C}[Th]{Corollary}

\author{Alexander Brudnyi\thanks{Research supported in part by NSERC.
}\\
Department of Mathematics and Statistics\\
University of Calgary, Calgary\\
Canada}
\title{Contractibility of Maximal Ideal Spaces of Certain Algebras of 
Almost Periodic Functions}
\date{} 
\begin{document} 
\maketitle
\begin{abstract}
{We study some topological properties of maximal ideal spaces of certain 
algebras of almost periodic functions. Our main result is that such spaces are
contractible. We present several analytic corollaries of this result.} 
\end{abstract}


\sect{\hspace*{-1em}. Introduction.}
{\bf 1.1.}
This paper is devoted to the study of some topological properties of maximal 
ideal spaces of certain algebras of almost periodic functions.
We apply our results to the solution of a matrix version of the corona 
problem for such algebras. To formulate these results we first introduce some 
notation and basic definitions.

The algebra $AP$ of (uniformly) {\em almost periodic functions} is, by
definition, the smallest closed subalgebra of $L^{\infty}(\Re)$ that contains 
all the functions $e_{\lambda}:=e^{i\lambda x}$, $\lambda\in\Re$.\\
A finite sum $p=\sum r_{j}e_{\lambda_{j}}\in AP$ with $r_{j}\in\Co$ and 
$\lambda_{j}\in\Re$ is called an {\em almost periodic polynomial}. 
The set of almost periodic polynomials is denoted by $AP^{0}$.

The {\em Bohr-Fourier spectrum} $\Omega(a)$ of
$a\in AP$ is the subset of $\Re$ defined by the formula
$$
\Omega(a):=\{\lambda\in\Re\ :\ 
\lim_{N\to\infty}\frac{1}{2N}\int_{-N}^{N}a(x)e_{-\lambda}(x)dx\neq 0\}\ .
$$
It is well known that $\Omega(a)$ is at most countable, see, e.g., 
[BKS, Sect. 2.7]. For instance, for 
$p=\sum_{j=1}^{k}r_{j}e_{\lambda_{j}}\in AP^{0}$ with all $r_{j}$ nonzero, 
$\Omega(p)=\{\lambda_{1},\dots,\lambda_{k}\}$.

Let $\Sigma\subset [0,\infty)$ be an
additive semigroup (i.e., $\lambda,\mu\in\Sigma$ implies that 
$\lambda+\mu\in\Sigma$) and suppose $0\in\Sigma$. Put 
$$
AP_{\Sigma}^{0}:=\{a\in AP^{0}\ :\ \Omega(a)\subset\Sigma\}\ .
$$
Let $AP_{\Sigma}$ be the closure of $AP_{\Sigma}^{0}$ in $AP$. Corollary
7.6 of [BKS] shows that
$$
AP_{\Sigma}=\{a\in AP\ :\ \Omega(a)\subset\Sigma\}\ .
$$
Since $\Sigma$ is a semigroup containing 0, $AP_{\Sigma}$ is a Banach 
algebra. Let $M(AP_{\Sigma})$ be the maximal ideal space of $AP_{\Sigma}$,
that is the space of nonzero homomorphisms \penalty-10000 $AP_{\Sigma}\to\Co$
equipped with the weak$*$ topology (i.e., the Gelfand topology).
The analytic structure of 
$M(AP_{\Sigma})$ is described by the Arens-Singer theorem, 
see, e.g, [BKS, Theorem 12.4]. In contrast, our main result describes the 
topological structure of $M(AP_{\Sigma})$.
\begin{Th}\label{te1}
$M(AP_{\Sigma})$ is contractible.
\end{Th}
\begin{C}\label{c1}
All \v{C}ech cohomology groups $H^{m}(M(AP_{\Sigma}),\Z)$, $m\geq 1$, 
are trivial.
\end{C}
{\bf 1.2.} Note that every function $a\in AP_{\Sigma}$ can be extended to a 
holomorphic function $\hat a$ in the upper half-plane $\H_{+}$ by the 
Poisson integral formula. Obviously, 
$$
\sup_{\H_{+}}|\hat a|=\sup_{\Re}|a|\ .
$$
This produces a natural embedding $\H_{+}\hookrightarrow M(AP_{\Sigma})$,
$z\mapsto\{${\em evaluation at} $z\}$.
The {\em corona problem} for $AP_{\Sigma}$ is to decide whether $\H_{+}$
is dense in $M(AP_{\Sigma})$ (in the Gelfand topology). The density 
is equivalent to the following statement (see, e.g., \penalty-10000
[BKS, Sect. 12.1]):

{\em Given a collection $f_{1},\dots, f_{n}$ of functions from 
$AP_{\Sigma}$ satisfying}
\begin{equation}\label{e1}
\inf_{z\in\H_{+}}\sum_{j=1}^{n}|\hat f_{j}(z)|>0\ ,
\end{equation}
{\em there are $g_{1},\dots,g_{n}\in AP_{\Sigma}$ such that}
$$
\hat f_{1}(z)\hat g_{1}(z)+\dots + \hat f_{n}(z)\hat g_{n}(z)=1\ \ \
{\rm for\ all}\ \ \ z\in\H_{+}\ .
$$
In [BKS, Theorem 12.7] a necessary and sufficient condition for
$\Sigma$ is given to conclude that $\H_{+}$ is dense in $M(AP_{\Sigma})$.
For instance, it is true if $\Sigma=\Delta\cap [0,\infty)$ where 
$\Delta\subset\Re$ is an additive group. (It is not true, e.g., for
$\Sigma=\{k+l\sqrt{2}\ :\ k,l\in\Z_{+}\}$, see \penalty-10000
[BKS, Example 12.9].)

In this paper we consider the following matrix version of the corona problem.

{\em Let $A=(a_{ij})$ be an $n\times k$ matrix, $k<n$, with entries in 
$AP_{\Sigma}$. Assume that
the family of determinants of the submatrices of $A$ of order $k$ satisfies
condition (\ref{e1}). Is there an $n\times n$ matrix
$\widetilde A=(\widetilde a_{ij})$, $\widetilde a_{ij}\in AP_{\Sigma}$,
so that $\widetilde a_{ij}=a_{ij}$ for $1\leq j\leq k$, $1\leq i\leq n$, and
$det(\widetilde A)=1$?}

In case that such a $\widetilde A$ exists we say that $\widetilde A$ completes
$A$. As a corollary of Theorem \ref{te1} we obtain the following 
result.\footnote{Theorem \ref{te2} answers a
question of L.Rodman [Ro].}
\begin{Th}\label{te2}
Suppose that $\H_{+}$ is dense in $M(AP_{\Sigma})$. Then for any $A$
satisfying the matrix version of the corona problem there is a matrix
$\widetilde A$ completing it.
\end{Th}

Finally, we formulate one other application of Theorem \ref{te1}.
\begin{Th}\label{te3}
Suppose that $\H_{+}$ is dense in $M(AP_{\Sigma})$. Let $f\in AP_{\Sigma}$ 
satisfy (\ref{e1}) (with $n=1$, $f_{1}:=f$).
Then there is $g\in AP_{\Sigma}$ such that $f=e^{g}$.
\end{Th}
\sect{\hspace*{-1em}. Proofs.}
{\bf Proof of Theorem \ref{te1}.}
Let $X(\Sigma)$ be the vector space of complex-valued functions
on $\Sigma$ equipped with the product topology. Consider the map
$j:\H_{+}\to X(\Sigma)$,
$$
j(z)(\lambda):=e^{i\lambda z}\ ,\ \ \ z\in\H_{+},\ \lambda\in\Sigma\ .
$$
Let $\{z_{\lambda}\}_{\lambda\in\Sigma}$ be the family of coordinate
functionals on $X(\Sigma)$, i.e., $z_{\lambda}(f):=f(\lambda)$, 
$f\in X(\Sigma)$. By $P(\Sigma)$ we denote the algebra
of polynomials on $X(\Sigma)$ in variables $\{z_{\lambda}\}$. By definition,
any $p\in P(\Sigma)$ is a finite sum 
$\sum r_{i_{1}\dots i_{k}}z_{\lambda_{i_{1}}}\cdots z_{\lambda_{i_{k}}}$ with 
$r_{i_{1}\dots i_{k}}\in\Co$. Let $H(\Sigma)\subset X(\Sigma)$ be the 
polynomial hull of $j(\H_{+})$ with respect to $P(\Sigma)$, that is, 
$$
z\in H(\Sigma)\ \ \ \Longleftrightarrow\ \ \
|p(z)|\leq\sup_{j(\H_{+})}|p|\ \ \ {\rm for\ any}\ \ \ 
p\in P(\Sigma)\ .
$$
Then the direct limit construction of Royden [R] implies that $H(\Sigma)$
is homeomorphic to $M(AP_{\Sigma})$. From now on we identify these two 
objects.

Further, define the map $R:X(\Sigma)\times [0,1]\to X(\Sigma)$ by the formula
$$
R(f,t)(\lambda):=t^{\lambda}f(\lambda)\ ,\ \ \ f\in X(\Sigma),\ t\in [0,1]\ ,
\ \lambda\in\Sigma\ .
$$
Clearly this map is continuous if we consider $X(\Sigma)$ with the 
product topology and $[0,1]$ with the usual one. So $R$ determines a
contraction of $X(\Sigma)$ to $0$. We will show that 
$R$ maps $M(AP_{\Sigma})\times [0,1]$ to $M(AP_{\Sigma})$. This will complete
the proof.

First, prove that $R$ maps $j(\H_{+})\times (0,1]$ to $j(\H_{+})$.
Indeed, for a fixed $t\in (0,1]$ we have
$$
R(j(z),t)(\lambda):=t^{\lambda}e^{i\lambda z}=e^{i\lambda(-i\log t)}
e^{i\lambda z}=e^{i\lambda (z+i\log(1/t))}=j(z+i\log(1/t))(\lambda)\ .
$$
We used here that $\log(1/t)\geq 0$. 
The above identity implies that $R(j(z),t)=j(z+i\log(1/t))$ as required. 
Note also that the closure $\overline{j(\H_{+})}$ contains $0$. \\
Now suppose that
there is $v\in M(AP_{\Sigma})$ and $t\in (0,1)$ such that 
$R(v,t)\not\in M(AP_{\Sigma})$. Then there is a polynomial 
$p(z_{\lambda_{1}},\dots, z_{\lambda_{k}})\in P(\Sigma)$ such that 
$$
|p(t^{\lambda_{1}}v_{\lambda_{1}},\dots, t^{\lambda_{k}}v_{\lambda_{k}})|>
\sup_{j(\H_{+})}|p|\
$$
where $v_{\lambda_{i}}:=z_{\lambda_{i}}(v)$.
Consider the polynomial $q(z_{\lambda_{1}},\dots, z_{\lambda_{k}}):=
p(t^{\lambda_{1}}z_{\lambda_{1}},\dots, t^{\lambda_{k}}z_{\lambda_{k}})$. 
Then the previous inequality implies that
$$
|q(v)|>\sup_{j(\H_{+})}|p|\geq\sup_{j(t\H_{+})}|p|=\sup_{j(\H_{+})}|q|\ .
$$
This contradicts to the fact that $v\in M(AP_{\Sigma})$.\ \ \ \ \ $\Box$\\
{\bf Proof of Corollary \ref{c1}.} The proof follows straightforwardly from
the fact that $M(AP_{\Sigma})$ is homotopic to a point.\ \ \ \ \ $\Box$\\
{\bf Proof of Theorem \ref{te2}.} Since $\H_{+}$ is dense in $M(AP_{\Sigma})$
and $M(AP_{\Sigma})$ is contractible, the proof is a direct consequence of
Theorem 3 of Lin [L].\ \ \ \ \ $\Box$\\
{\bf Proof of Theorem \ref{te3}.} The proof follows directly from the
Arens-Royden theorem (see [R]), because in our case 
$H^{1}(M(AP_{\Sigma}),\Z)=0$.\ \ \ \ \ $\Box$
%===================


\begin{thebibliography}{}
\bibitem[BKS]{BKS}
A. B\"{o}ttcher, Yu. Karlovich and I. Spitkovsky, Convolution operators
and factorization of almost periodic matrix functions, Operator Theory:
Advances and Applications, {\bf 131}, Birkh\"{a}user Verlag, Basel, 2002.
\bibitem[L]{L}
V. Lin, Holomorphic fibering and multivalued functions of elements of a
Banach algebra, Func. Anal. Appl., {\bf 7} (2) (1973), 122-128, English 
translation.
\bibitem[R]{R}
H. Royden, Function algebras, Bull. Amer. Math. Soc., {\bf 69} (1963),
281-298. 
\bibitem[Ro]{Ro}
L. Rodman, Private communication.

\end{thebibliography}
\end{document}